\documentclass[11pt]{article}

\usepackage{url}
\usepackage{doi}

\hypersetup{pdfproducer={other system},pdfcreator={other tool}}
\hypersetup{%
linkcolor=blue,anchorcolor=blue,%
citecolor=blue,filecolor=blue,%
menucolor=blue,%
urlcolor=blue,%
colorlinks=true}

\begin{document}



\begin{center} \Large
Recent implementations, applications, and extensions of the 
Locally Optimal Block Preconditioned Conjugate Gradient method (LOBPCG)
\end{center}

\begin{center} \em
\underline{Andrew Knyazev (Mitsubishi Electric Research Laboratories)}
\end{center}


\centerline{Abstract}
\medskip
Since introduction [A. Knyazev, Toward the optimal preconditioned eigensolver: Locally optimal block preconditioned conjugate gradient method, SISC (2001) \doi{10.1137/S1064827500366124}] and efficient parallel implementation [A. Knyazev et al., Block locally optimal preconditioned eigenvalue xolvers (BLOPEX) in HYPRE and PETSc, SISC (2007) \doi{10.1137/060661624}], LOBPCG has been used is a wide range of applications in mechanics, material sciences, and data sciences. We review its recent implementations and applications, as well as extensions of the local optimality idea beyond standard eigenvalue problems. 

\section{Background}
Kantorovich in 1948 has proposed calculating the smallest eigenvalue $\lambda_1$ of a symmetric matrix $A$ 
by steepest descent using a direction $r=Ax-\lambda(x)$ of a scaled gradient of a 
Rayleigh quotient $\lambda (x) = (x,A x) / (x, x)$ in a scalar product $(x,y)=x^\prime y$, 
where the step size is computed by minimizing the Rayleigh quotient in the span of the vectors $x$ and $w$, i.e. in a 
locally optimal manner. Samokish in 1958 has proposed applying a preconditioner $T$ to the vector $r$ to generate the preconditioned direction $w=Tr$ and derived asymptotic, as $x$ approaches the eigenvector, convergence rate bounds.
Block locally optimal multi-step steepest descent is described in Cullum and Willoughby in 1985.

Local minimization of the Rayleigh quotient on the subspace spanned by the current approximation, the current residual and the previous approximation, as well as its block version, appear in AK PhD thesis; see {1986}. The preconditioned version is analyzed in 1991 and 1998, and a ``practically stable'' implementation \cite{K01}. 
LOBPCG general technology can also be viewed as a particular case of generalized block Davidson diagonalization methods with thick restart, or accelerated block gradient descent with plane-search. 
The main points on LOBPCG are as follows: 
 \begin{itemize}
 \item The costs per iteration and the memory use in LOBPCG are competitive with those of the Lanczos method, computing a single extreme eigenpair. 
 \item Linear convergence is theoretically guaranteed and practically observed, since
local optimality  implies that LOBPCG converges at least as fast as the gradient descent. In numerical tests,  LOBPCG typically shows no super-linear convergence.
 \item LOBPCG blocking allows utilizing highly efficient matrix-matrix operations, e.g., BLAS 3.
 \item LOBPCG can directly take advantage of preconditioning, in contrast to the Lanczos method.
 LOBPCG allows variable and non-symmetric or positive definite preconditioning.   
 \item LOBPCG allows warm starts and computes an approximation to the eigenvector on every iteration. It has no numerical stability issues similar to those of the Lanczos method.
 \item LOBPCG is reasonably easy to implement, so many implementations have appeared.
 \item Very large block sizes in LOBPCG become expensive to deal with due to orthogonalizations and the use of the Rayleigh-Ritz method on every iteration.
\end{itemize}

\section{General purpose open source public software implementations}
 \begin{itemize}
 \item MATLAB/Octave version by AK of the reference algorithm from 2001, publicly available 
 from MATLAB/Octave repositories. 
\item BLOPEX, developed by AK and his team \cite{BLOPEX}. An abstract implementation of the reference algorithm from 2001.
Originally at Google Code, currently hosted at bitbucket.org by Jose E. Roman.
\item BLOPEX included into SLEPc/PETSc by Jose E. Roman.
\item BLOPEX incorporated into hypre by AK and his team with hypre developers since 2007.
\item LOBPCG with orthogonalizations, Anasazi/Trilinos by U. Hetmaniuk and R. Lehoucq since 2006.
\item LOBPCG in Java at Google Code Archive by M. E. Argentati. 
\item LOBPCG in Python: SciPy, scikit-learn and megaman for manifold learning.
\item LOBPCG GPU implementation in NVIDIA AmgX and spectral graph partitioning.
\item LOBPCG heterogeneous CPU-GPU implementation by H. Anzt, S. Tomov, and J. Dongarra using blocked sparse matrix vector products in MAGMA (2015).
\end{itemize}

\section{LOBPCG in applications}
As a generic eigenvalue solver for large-scale symmetric eigenproblems, LOBPCG is being used in a wide range
of applications. Our own small contribution is to data clustering and image segmentation (2003) \cite{k2003} and (2015) 
\cite{20150363361}, as well as for low-pass graph-based signal filtering (2015) \cite{KM15a}. Below we list some application-oriented open source software implementations of LOBPCG in various domains.

 Material Sciences:
 \begin{itemize}
 \item ABINIT (including CUDA version) by AK with ABINIT developers, since 2008. Implements density functional theory, using a plane wave basis set and pseudopotentials.
 \item Octopus TDDFT, employs pseudopotentials and real-space numerical grids.
 \item Pescan, nonselfconsistent calculations of electron/hole states in a plane wave basis set.
 \item Gordon Bell Prize finalist at ACM/IEEE Conferences on Supercomputing in 2005, \cite{1559996}.
  \item Gordon Bell Prize finalist at ACM/IEEE Conferences on Supercomputing in 2006, \cite{Yamada:2006:HCE:1188455.1188504}.
  \item KSSOLV A MATLAB Toolbox for Solving the Kohn Sham Equations, Yang (2009).
 \item TTPY, spectra of molecules using tensor train decomposition, Rakhuba, Oseledets (2016).
 \item Platypus‐QM, PLATform for dYnamic protein unified simulation, Takano et al. (2016).
 \item MFDn,  Nuclear Configuration Interaction, Shao et al. (2016)
\end{itemize}

Data mining:
 \begin{itemize}
 \item sklearn: Machine Learning in Python, spectral clustering.  
 \item Megaman: Scalable Manifold Learning in Python, McQueen et al. (2016).
 \end{itemize}

Multi-physics:
 \begin{itemize}
 \item NGSolve (Joachim Schoberl), a general purpose adaptive finite element library with multigrid.
 \end{itemize}

Electromagnetic calculations:
  \begin{itemize}
 \item NGSolve by Joachim Schoberl Maxwell solvers with Python interface.
 \item MFEM by Tzanio Kolev, highly scalable parallel high-order adaptive finite elements interfaced to hypre LOBPCG and multigrid  preconditioning.
 \item Python based PYFEMax by R. Geus, P. Arbenz , L. Stingelin. Nedelec finite element discretisation of Maxwell’s equations.
 \end{itemize}

\section{Methods, motivated by LOBPCG}
Suppose that the number $n_b$ of vectors in the block must be small, e.g., due to memory limitations. 
Let us consider the extreme case $n_b=1$ of single-vector iterations only.
At the same time, assume that we need to compute $n_v>n_b$ eigenpairs.
It is well-known that one can compute $n_b$ eigenpairs at the time, lock them together with all the previously 
computed eigenvectors, and sequentially compute the next $n_b$ eigenpairs in the orthogonal complement to all already locked approximate eigenvectors. E.g., if $n_b=1$, the eigenpairs are computed one-by-one sequentially. 
In contrast, LOBPCG II, see \cite{K01} utilities $n_v$ parallel execution of single vector ($n_b=1$) 3-term locally optimal recurrences combined with the Rayleigh-Ritz procedure on the $n_v$ -dimensional subspace spanned by the current approximations 
to all $n_v$ eigenvectors. E. Vecharynski, C. Yang, and J. Pask (2015) extend the LOBPCG II approach to $n_v/n_b$ parallel execution of $n_b>1$ block-vector 3-term recurrences with additional orthogonalizations, and demonstrate that Rayleigh-Ritz can be omitted on many iterations, in Quantum Espresso plane-wave DFT electronic structure software.

There are various other recent developments, related to LOBPCG: 
 \begin{itemize}
\item Indefinite variant of LOBPCG for definite matrix pencils, D. Kressner (2013) is based on an indefinite inner product, and generalizes LOBP4DCG method by Bai and Li (2012 and 2013) for solving product eigenvalue problems
\item {Folded spectrum to LOBPCG (FS-LOBPCG), V.S. Sunderam (2005)} adds three more block vectors that store the matrix-vector products of the blocks $X$, $R$, and $P$.
\item Szyld and Xue (2016) extend LOBPCG to nonlinear Hermitian eigenproblems with variational characterizations.
\end{itemize}

{Low-rank tensor related variants of LOBPCG:}

 \begin{itemize}
\item Hierarchical Tucker decomposition in low-rank LOBPCG, D. Kressner at al., 2011, 2016.
\item Tensor-Train (TT) manifold-preconditioned LOBPCG for many-body Schr\"{o}dinger equation, Oseledets at al., 2014, 2016
\end{itemize}

\section{Conclusions}
 Still nobody can pronounce it, but LOBPCG:
 \begin{itemize}
 \item has now a life on its own, with ~7K hits in google search,
 \item accepted as one of standard methods in DFT electronic structure calculations,
 \item appears in most major HPC open source libraries: PETSc, hypre and  Anasazi/Trilinos, e.g., pre-packaged in OpenHPC cluster management software
 \item implemented for GPUs for ABINIT and in MAGMA and AmgX
 \end{itemize}

 More efforts are needed to decrease the use of orthogonalizations and Rayleigh-Ritz for large block sizes, while still preserving fast convergence, maintaining reasonable stability, and keeping the algorithm simple to implement.


\end{document}